\documentclass[a4paper]{article}     % `article' with A4 paper size option
\usepackage{amsthm}                 % theorem extensions
\usepackage{amssymb}
\usepackage{amsmath}
\usepackage{ae}
\usepackage[latin1]{inputenc}
\usepackage[T1]{fontenc}
\usepackage[all]{xy}
\newcommand{\sign}{\textnormal{sign}}
\newcommand{\ind}{\textnormal{ind}}
\newcommand{\rg}{\textnormal{rank}}

\newcommand{\res}{\textnormal{res}}

\newcommand{\mo}{\textnormal{mod}}

\newcommand{\ann}{\textnormal{ann}}

\newcommand{\soc}{\textnormal{soc}}

\newtheorem{theorem}{Theorem}
\newtheorem*{theorem*}{Theorem}
\newtheorem{corollary}{Corollary}
\newtheorem{proposition}{Proposition}
\newtheorem{lemma}{Lemma}

\newtheorem{claim}{Claim}
\newcommand{\CC}{{\mathbb{C}}}

\newcommand{\RR}{{\mathbb{R}}}
\newcommand{\ZZ}{{\mathbb{Z}}}

\newcommand{\calO}{{\mathcal{O}}}
\newcommand{\calB}{{\mathcal{B}}}

\newcommand{\calA}{{\cal A}}
\newcommand{\calC}{{\cal C}}

\newcommand{\OCn}{\mathcal{O}_{\CC^n,0}}

\newcommand{\ICn}{\ind_{\CC^n,0}}
\newcommand{\IVn}{\ind_{V,0}}

\newcommand{\RCn}{\res_{\CC^n,0}}

\newcommand{\RVn}{\res_{V,0}}
\newcommand{\RCpi}{\res_{\CC^n,p_i}}

\begin{document}

\title{Local residues of holomorphic 1-forms on an isolated surface singularity}
\author{Oliver Klehn\\
\parbox{9cm}{\small
 \begin{center} Institut f\"{u}r Mathematik, Universit\"{a}t Hannover, \\
        Postfach 6009, D-30060 Hannover, Germany \\
        E-mail: klehn@math.uni-hannover.de
\end{center}}
\date{}
}

\maketitle

\begin{abstract}
We define local residues of holomorphic 1-forms on an isolated surface singularity that have isolated zeros and
prove that a certain residue equals the index of the 1-forms defined by Ebeling and Gusein-Zade.
\end{abstract}

\section{Introduction}
\label{intro}
Let $f_1,\dots ,f_n\in\OCn$ define a regular sequence, let $J_f$ be the determinant of the Jacobian matrix of
$f:=(f_1,\dots ,f_n)$ and set
\[
Q_f:=\frac{\OCn}{(f_1,\dots ,f_n)}.
\]
Then a classical result states
\[
\ICn f=\dim_{\CC}Q_f=\RCn
\begin{bmatrix}
J_f\\
f_1 \dots f_n
\end{bmatrix},
\]
where $\ind $ is the Poincar\'{e}-Hopf index and $\res $ the local residue symbol. Recall that the residue is defined
for any $h\in\OCn$ by
\[
\RCn \begin{bmatrix} h\\ f_1\dots f_n \end{bmatrix}:=\frac{1}{(2\pi i)^n}\int_{\Gamma}
\frac{hdz_1\wedge\dots\wedge dz_n}{f_1\dots f_n}
\]
where $\Gamma$ is the real $n$-cycle $\{ |f_i|=\epsilon_i ,i=1,\dots ,n\}$ for $\epsilon_i\in\RR_{>0}$ chosen small enough
oriented so that $d( \arg f_1)\wedge\dots\wedge d(\arg f_n)\geq 0$, and the Poincar\'{e}-Hopf index is the degree
of the map $f/||f||$ restricted to a small sphere at the origin, see also \cite{agv,gh,m} for the definitions.

In the framework of singularity theory it
arises a natural question: How to generalize these definitions when one considers holomorphic function germs on an
isolated singularity $(V,0)\subset (\CC^n,0)$? There are two directions one can go. The first is the definition
of an index of holomorphic vector fields tangent to $(V,0)$, called the GSV index. The results can be found in
\cite{bg,g,gsv,ls1,ls2,lss,ss}. The other direction is to consider holomorphic 1-forms on $(V,0)$ which has been done
by Ebeling and Gusein-Zade in \cite{eg}. This is the situation we give a generalization of the residue symbol
and therefore we briefly recall the main results of \cite{eg}.

Let $(V,0)=(\{ f_1=\dots =f_q=0\},0)\subset (\CC^n,0)$ be an isolated complete intersection singularity and
$\omega =\sum_{i=1}^n\omega_idz_i$ the germ of a holomorphic 1-form on $(\CC^n,0)$ which has an isolated zero on $(V,0)$.
Choose a sufficiently small sphere $S_{\delta}$ around the origin in $\CC^n$ which intersects $V$ transversally
and consider the link $K=V\cap S_{\delta}$ of $V$. The 1-forms $\omega,df_1,\dots ,df_q$ are linearly independent for all points
of $K$ and we have a well defined map
\[
(\omega, df_1,\dots ,df_q)\colon K\to W_{q+1}(\CC^n)
\]
where $W_{q+1}(\CC^n)$ denotes the manifold of $(q+1)$-frames in the dual $\CC^n$. We have
\[
H_{2n-2q-1}(K)\cong\ZZ,\;\;H_{2n-2q-1}(W_{q+1}(\CC^n))\cong\ZZ
\]
and therefore the map has a degree. We let $K$ to be oriented as boundary of the complex manifold $V\setminus\{ 0\}$ here.
The index $\IVn\omega$ of $\omega$ is defined to be the degree of this map. (If $V$ is a curve $K$ can have more components, we will sum
over the degrees of the components then.)

Let $J$ be the ideal in $\OCn$ generated by $f_1,\dots ,f_q$ and the $(q+1)$-minors of the matrix
\[
\left(\begin{array}{ccc}
\frac{\partial f_1}{\partial z_1} & \dots &\frac{\partial f_1}{\partial z_n}\\
\vdots                                      & \ddots& \vdots                           \\
\frac{\partial f_q}{\partial z_1}  & \dots &\frac{\partial f_q}{\partial z_n}\\
\omega_1                                     & \dots & \omega_n
\end{array}\right).
\]
For a regular value $0\neq\epsilon=(\epsilon_1,\dots ,\epsilon_q)\in\CC^q$ of $f$ chosen small enough and
a small ball $B_{\delta}$ around the origin in $\CC^n$ define
$V_{\epsilon}:=f^{-1}(\epsilon )\cap B_{\delta}$. $V_{\epsilon}$ is transversally to $S_{\delta}$ then. We set $\mathcal{A}:=\OCn /J$
and have the following theorem of Ebeling and Gusein-Zade:

\begin{theorem}
(i) If $(V,0)$ is smooth, $\IVn \omega$ is the usual Poincar\'{e}-Hopf index.\\
(ii) $\IVn \omega$ equals the number of zeros of $\omega$ on $V_{\epsilon}$ counted with multiplicities.\\
(iii) $\IVn \omega =\dim_{\CC}\mathcal{A}$.
\end{theorem}

In this paper we define if $V$ is a surface a linear form
\[
\RVn \colon \calA \to \CC,
\]
which we call the relative residue form of $\omega$ and prove that for a certain class $\sigma\in\calA$ one has a formula
\[
\IVn \omega =\RVn (\sigma ).
\]
We use the linear form to prove some algebraic properties of $\calA$. There is used classical Grothendieck residue
theory and we will therefore recall the main facts as them can be found in \cite{gh}. With the notations introduced
at the beginning of this section we summarize some results of \cite{gh} in the following theorem:

\begin{theorem}
(i) For $i=1,\dots ,n$ let $g_i:=\sum_{j=1}^n a_{ij}f_j$ so that $g_1,\dots ,g_n$ is a regular sequence and let $A$ be the matrix
$(a_{ij})_{i=1,\dots ,n}^{j=1,\dots ,n}$. Then one has for any $h\in\OCn$
\[
\RCn
\begin{bmatrix}
h\\ f_1 \dots f_n
\end{bmatrix} = \RCn
\begin{bmatrix}
h \det A\\ g_1 \dots g_n
\end{bmatrix}.
\]
(ii) The residue defines a linear form
\[
\RCn
\begin{bmatrix}
\cdot \\ f_1\dots f_n
\end{bmatrix}
\colon Q_f \to\CC .
\]
(iii)
The induced pairing
\[
B\colon Q_f\times Q_f\to\CC
\]
defined by
\[
B(h,g):=\RCn
\begin{bmatrix}
h\cdot g\\ f_1\dots f_n
\end{bmatrix}
\]
is non-degenerate.\\
(iv)
For perturbations $f_{\epsilon}$ and $h_{\epsilon}$ of $f$ and $h$ one has
\[
\RCn
\begin{bmatrix} h\\ f_1 \dots f_n \end{bmatrix} = \lim_{\epsilon\to 0} \sum_i \RCpi
\begin{bmatrix} h_{\epsilon} \\ f_{\epsilon,1} \dots f_{\epsilon,n} \end{bmatrix}
\]
where one sums over the zeros of $f_{\epsilon}$ in a neighbourhood of the origin.
\end{theorem}

We want to formulate our main result. Let $q=n-2$ and $M_i$ the matrix obtained from
\[
\left(\begin{array}{ccc}
\frac{\partial f_1}{\partial z_1} & \dots &\frac{\partial f_1}{\partial z_n}\\
\vdots                                      & \ddots& \vdots                           \\
\frac{\partial f_{n-2}}{\partial z_1}  & \dots &\frac{\partial f_{n-2}}{\partial z_n}\\
\omega_1                                     & \dots & \omega_n
\end{array}\right)
\]
by cancelling the $i$-th column, $m_i:=\det (M_i)$ and $M:=\left( (-1)^i\frac{\partial m_i}{\partial z_j}\right )$.
Let $\sigma$ be the coefficient of $t^{n-2}$ in the characteristical polynomial of $M$. We will prove that there
exits a linear change of coordinates (in fact a generic one) so that $(m_1,m_2)$ is a regular $\mathcal{O}_{V,0}$-sequence
and call these coordinates good coordinates. Moreover we define for $h\in\OCn$
\[
\RVn \begin{bmatrix} h\\ m_1 m_2 \end{bmatrix} :=\frac{1}{(2\pi i)^2}\int_{\Sigma}
\frac{ hdz_1\wedge dz_2}{m_1 m_2}
\]
with $\Sigma :=\{f_1=\dots =f_{n-2}=0, |m_1|=\delta_1 ,|m_2|=\delta_2 \}$ oriented so that $d(\arg m_1)\wedge d(\arg m_2)\geq 0$
and where $\delta_1 ,\delta_2$ are small positive real numbers. The following theorem is our main result:

\begin{theorem}
\label{thm1}
In each system of good coordinates we have a linear form
\[
\RVn \begin{bmatrix}\cdot \\ m_1 m_2\end{bmatrix}\colon\calA\to\CC
\]
with
\[
\dim_{\CC}\calA=\RVn
\begin{bmatrix}
\sigma\\  m_1 m_2
\end{bmatrix}.
\]
\end{theorem}

I wish to thank W. Ebeling and S.M. Gusein-Zade for useful discussions and the referee for useful comments
concerning the presentation. The idea that one should try to find a residue formula
for the index of a holomorphic 1-form is due to W. Ebeling.

\section{Absolute and relative residues}
\label{sec1}

Let $(V,0)\subset (\CC^n,0)$ be an icis as before and $g_1,\dots ,g_{n-q}\in\OCn$ define an isolated zero
on $(V,0)$. If the real hypersurfaces $\{|g_i|=\delta_i\}$ for small $\delta_i$ are in general position let
\[
\Sigma :=\{f_1=\dots =f_q=0, |g_i|:=\delta_i ,i=1,\dots ,n-q\}
\]
to be the real $(n-q)$-cycle oriented so that $d(\arg g_1)\wedge\dots \wedge d(\arg g_{n-q})\geq 0$.
Then we have for each
\[
\eta\in\Omega_{V,0}^{n-q}:=\Omega_{\CC^n,0}^{n-q}/
(\Sigma f_i\Omega_{\CC^n,0}^{n-q}+\Sigma df_i\wedge\Omega_{\CC^n,0}^{n-q-1})
\]
a well defined integral
\[
\frac{1}{(2\pi i)^{n-q}}\int_{\Sigma}\frac{\eta}{g_1\dots g_{n-q}}
\]
which we denote by
\[
\RVn
\begin{bmatrix}
\eta \\ g_1\dots g_{n-q}
\end{bmatrix}.
\]
We now want to prove a relation between these relative residues and the absolute residues. Define
\[
\lambda\colon\Omega_{V,0}^{n-q}\to\calO_{V,0}
\]
by $\lambda (\eta):=h$, where $h$ is defined by
\[
hdz_1\wedge\dots\wedge dz_n:=\eta\wedge df_1\wedge\dots\wedge df_q.
\]
We want to prove the following theorem:

\begin{theorem}
\label{thm2}
\[
\RVn
\begin{bmatrix}
\eta \\ g_1\dots g_{n-q}
\end{bmatrix} = \RCn
\begin{bmatrix}
\lambda (\eta) \\ g_1\dots g_{n-q} f_1\dots f_q
\end{bmatrix}.
\]
\end{theorem}

The intersection multiplicity $I$ of the hypersurfaces defined by $f_1,\dots ,f_q$, $g_1,\dots ,g_{n-q}$ is also given as a relative integral:

\begin{corollary}
\label{cormult}
\[
I= \frac{1}{(2\pi i)^{n-q}}\int_{\Sigma}\frac{dg_1\wedge\dots\wedge dg_{n-q}}{g_1\dots g_{n-q}}.
\]
\end{corollary}

The Theorem follows from the next Lemma if we set
\[
DF:=\det\left( \frac{\partial (f_1,\dots ,f_q)}{\partial (z_{n-q+1},\dots ,z_n)}\right) .
\]

\begin{lemma}
\label{res}
\[
\RVn
\begin{bmatrix}
hdz_1\wedge\dots\wedge dz_{n-q} \\ g_1\dots g_{n-q}
\end{bmatrix} = \RCn
\begin{bmatrix}
h\cdot DF \\ g_1\dots g_{n-q} f_1\dots f_q
\end{bmatrix}.
\]
\end{lemma}

\begin{proof}[Proof of Theorem \ref{thm2}]
Let $1\leq i_1<\dots <i_{n-q}\leq n$ and $1\leq j_1<\dots <j_q\leq n$ the complement,
$\sigma :=(i_1,\dots ,i_{n-q},j_1,\dots ,j_q)\in S_n$.
Then we get
\begin{equation}
\label{eqres}
\RVn
\begin{bmatrix}
hdz_{i_1}\wedge\dots\wedge dz_{i_{n-q}} \\ g_1\dots g_{n-q}
\end{bmatrix} = \sign\sigma\cdot\RCn
\begin{bmatrix}
h\cdot
\frac{\partial (f_1,\dots ,f_q)}{\partial (z_{j_1},\dots ,z_{j_q})} \\ g_1\dots g_{n-q} f_1\dots f_q
\end{bmatrix}
\end{equation}
by Lemma \ref{res} using the permutation of coordinates $\sigma$ and the integral transformation formula. Now the theorem follows
from equation \ref{eqres} using the Laplace expansion formula.
\end{proof}

First we prove a special case of Lemma \ref{res}. Define
\[
F:=(g_1,\dots ,g_{n-q},f_1,\dots ,f_q).
\]
\begin{claim}
Lemma \ref{res} holds if $0$ is a regular value of $F$.
\end{claim}

To prove this we need some facts from linear algebra.
\begin{lemma}
\label{det1}
Let
$H:=
\left(\begin{array}{cc}
A & B\\
C & D
\end{array}\right)
$ be the decomposition of a $n\times n$-matrix in four blocks where $A$ and $D$ are squared and $A$ is invertible. Then
\[
\det H= \det A\cdot \det (D-CA^{-1}B).
\]
\end{lemma}

\begin{proof}
It is just a simple exercise.
\end{proof}

\begin{lemma}
\label{det2}
Let
$H:=
\left(\begin{array}{cc}
A & B\\
C & D
\end{array}\right)
$ be an invertible $n\times n$-matrix with $H^{-1}=
\left(\begin{array}{cc}
E & F\\
G & I
\end{array}\right)$ where $A$ und $E$ are $j\times j$-matrices and $D$ and $I$ are $(n-j)\times (n-j)$-matrices. Then
\[
\det D=\det E\cdot\det H.
\]
\end{lemma}

\begin{proof}
We have
\[
1=HH^{-1}=
\left(\begin{array}{cc}
AE+BG & AF+BI\\
CE+DG & CF+DI
\end{array}\right).
\]
Let $E$ be invertible. Sice we have $CE+DG=0$ we get $-CF-DGE^{-1}F=0$ and therefore with $DI=1-CF$ it follows
\begin{equation}
\label{gl}
\det (1-CF-DGE^{-1}F)=\det D\cdot\det (I-GE^{-1}F)=1.
\end{equation}
This means that $D$ is invertible. Now assume that $D$ is invertible. Since we have $CE+DG=0$ we get $BD^{-1}CE+BG=0$. From $1-BG=AE$ it follows that
\[
\det (1-BD^{-1}CE-BG)=\det (A-BD^{-1}C)\cdot\det E=1
\]
and so we find that $E$ must be invertible. We have shown that $E$ is invertible if and only if $D$ is invertible and so the Lemma follows
for non invertible $E$. If $E$ is invertible we get by Lemma  \ref{det1}
\[
\det H^{-1}=\det E\cdot\det (I-GE^{-1}F).
\]
By equation \ref{gl} we get
\[
\det H^{-1}=\det E\cdot\frac{1}{\det D}
\]
and $\det D=\det E\cdot\det H$ follows.
\end{proof}

\begin{proof}[Proof of the Claim]
We find that $V$ is smooth and $0$ is also a regular value of $G:=(g_1,\dots ,g_{n-q})|_V$. We choose new coordinates
$x:=F(z)$ on $(\CC^n,0)$ resp. $x':=G(z)$ on $(V,0)$ and set $\tilde{F}:=F^{-1}$ resp. $\tilde{G}:=G^{-1}$. Using the integral transformation formula
and the Cauchy integral formula we get
\small{
\begin{equation*}
\begin{split}
\RVn
\begin{bmatrix}
hdz_1\wedge\dots\wedge dz_{n-q} \\ g_1\dots g_{n-q}
\end{bmatrix} &=
\frac{1}{(2\pi i)^{n-q}}\int_{\Sigma}\frac{hdz_1\wedge\dots\wedge dz_{n-q}}{g_1\dots g_{n-q}}\\
&= \frac{1}{(2\pi i)^{n-q}}\int_T\frac{h(\tilde{G})\frac{\partial (\tilde{G}_1,\dots ,\tilde{G}_{n-q})}
{\partial (x_1,\dots ,x_{n-q})}dx_1\wedge\dots\wedge dx_{n-q}}{x_1\dots x_{n-q}}\\
&= h(\tilde{G}(0))\frac{\partial (\tilde{G}_1,\dots ,\tilde{G}_{n-q})}
{\partial (x_1,\dots ,x_{n-q})}(0)\\
&=   h(0)\frac{\partial (\tilde{F}_1,\dots ,\tilde{F}_{n-q})}
{\partial (x_1,\dots ,x_{n-q})}(0).
\end{split}
\end{equation*}
}
\normalsize{
$T$ is the torus
\[
T:=\{|x_i|=\delta_i,\; i=1,\dots ,n-q\}
\]
oriented as usually here. Similarly we get
\[
\RCn
\begin{bmatrix} h\cdot DF \\ g_1\dots g_{n-q}f_1\dots f_q \end{bmatrix}=
h(0)DF(0)\frac{\partial (\tilde{F}_1,\dots ,\tilde{F}_n)}
{\partial (x_1,\dots ,x_n)}(0).
\]
Now the equality of the residues follows from Lemma \ref{det2}.
 }\end{proof}

\begin{proof}[Proof of Lemma \ref{res}]
Denote by $\mu$ the local multiplicity of $F$. Let $G$ be defined as above. The regular values of $G$ are dense and we have for those
$y=(y_1,\dots ,y_{n-q})$ with $g_y:=(g_1-y_1,\dots ,g_{n-q}-y_{n-q})$
\[
\RCn \begin{bmatrix} h\cdot DF\\ g_1\dots g_{n-q} f_1\dots f_q \end{bmatrix}
=\lim_{y\to 0}\sum_{i=1}^{\mu}\RCpi
\begin{bmatrix} h\cdot DF\\ g_{y,1}\dots g_{y,n-q}f_1\dots f_q \end{bmatrix}.
\]
The value $(y,0)$ is a regular value of $F$: The preimages of $y$ under $G$ are $\mu$ simple zeros of $G-y$, which
are also the preimages of $(y,0)$ under $F$. They must be also simple zeros of $F-(y,0)$ then. By the Claim we get
\[
\RCn \begin{bmatrix} h\cdot DF\\ g_1\dots g_{n-q} f_1\dots f_q \end{bmatrix}
=\lim_{y\to 0}\frac{1}{(2\pi i)^{n-q}}\int_{\Sigma_{y,\rho}}
\frac{hdz_1\wedge\dots\wedge dz_{n-q}}{g_{y,1}\dots g_{y,n-q}}
\]
where
\[
\Sigma_{y,\rho}:=\{f_1=\dots =f_q=0,\;|g_{y,i}|=\rho_i ,i=1,\dots ,n-q\}
\]
and $\rho$ is chosen small enough. $\Sigma_{y,\rho}$ decomposes into $\mu$ components then. Now it sufficies to show that $\Sigma_{y,\rho}$
is homologous to $\Sigma$ in
\[
V\setminus \Bigl( \bigcup_{i=1}^{n-q}\{ g_i-y_i=0\}\cup \{ 0\}\Bigr).
\]
This can be done if $y$ is chosen small enough in the same lines as in \cite[p. 113]{agv}.
\end{proof}

\section{Residues of holomorphic 1-forms}
\label{sec2}
To prove our results on holomorphic 1-forms on an isolated surface singularity we
will look at the behaviour of the 1-forms on the Milnor fibre. Therefore we first study the smooth case. Then
we construct a class $\sigma\in\calA$ and show that it plays a similar role the Jacobian plays in the classical
case in the algebra $Q_f$. The Jacobian generates the $1$-dimensional socle of $Q_f$ and thus every linear form
$l\colon Q_f\to\CC$ which maps the Jacobian not to $0$ induces a non-degenerate pairing on $Q_f$
(The induced pairing is always the pairing defined by $B(h,g):=l(h\cdot g)$ in this paper). We will define a linear form
$\res\colon\calA\to\CC$ and use this form to prove a result on the dimension of the socle of $\calA$.

\subsection{The smooth case}
\label{sec21}
Let $(V,0)\subset (\CC^n,0)$ be an icis as before and $\omega=\sum_{i=1}^n\omega_idz_i$ the germ of a holomorphic 1-form on
$(\CC^n,0)$ with an isolated zero on $(V,0)$ (does not vanish on the tangent spaces $T_pV$ for $p\neq 0$ in a
neighbourhood of the origin). For $1\leq j_1,\dots ,j_{q+1}\leq n$ we set
\[
m_{j_1,\dots ,j_{q+1}}:=
\left|
\begin{array}{ccc}
\frac{\partial f_1}{\partial z_{j_1}} & \dots &\frac{\partial f_1}{\partial z_{j_{q+1}}}\\
\vdots                                      & \ddots& \vdots                                        \\
\frac{\partial f_q}{\partial z_{j_1}}  & \dots &\frac{\partial f_q}{\partial z_{j_{q+1}}}\\
\omega_{j_1}                                     & \dots & \omega_{j_{q+1}}
\end{array}
\right|.
\]
Furthermore for $1\leq i\leq n-q$ we set $\tilde{m}_i:=m_{i,n-q+1,\dots ,n}$. Let $J$ be the ideal as in the introduction and
$I$ the ideal in
$\OCn$ generated by $f_1,\dots ,f_q$ and the minors $\tilde{m}_1,\dots ,\tilde{m}_{n-q}$. For simplicity we may assume that
all minors vanish at the origin (otherwise the index is $0$).

\begin{lemma}
\label{eq1}
For any $i_1,\dots ,i_q,j_1,\dots ,j_{q+1}\in\{1,\dots ,n\}$ we have
\[
\frac{\partial (f_1,\dots ,f_q)}{\partial (z_{i_1},\dots ,z_{i_q})} m_{j_1,\dots ,j_{q+1}}=
\sum_{l=1}^{q+1} (-1)^{l+1}\frac{\partial (f_1,\dots ,f_q)}{\partial (z_{j_1},\dots ,\hat{z_{j_l}},\dots ,z_{j_{q+1}})} m_{j_l,i_1,\dots ,i_q}.
\]
\end{lemma}

\begin{proof}
Expansion of  $m_{j_1,\dots ,j_{q+1}}$ and $m_{j_l,i_1,\dots ,i_q}$ by the last row gives us
\begin{eqnarray*}
\frac{\partial (f_1,\dots ,f_q)}{\partial (z_{i_1},\dots ,z_{i_q})}  m_{j_1,\dots ,j_{q+1}} &=&
\sum_{l=1}^{q+1} (-1)^{l+1} \frac{\partial (f_1,\dots ,f_q)}{\partial (z_{j_1},\dots ,\hat{z_{j_l}},\dots ,z_{j_{q+1}})}  \Bigl(
m_{j_l,i_1,\dots ,i_q}\\
& - &\sum_{k=1}^q(-1)^{q+k}\omega_{i_k}
\frac{\partial (f_1,\dots ,f_q)}{\partial (z_{j_l},z_{i_1},\dots ,\hat{z_{i_k}},\dots ,z_{i_q})}\Bigr).
\end{eqnarray*}
Now it sufficies to show for any fixed $k$
\[
\sum_{l=1}^{q+1} (-1)^{l} \frac{\partial (f_1,\dots ,f_q)}{\partial (z_{j_1},\dots ,\hat{z_{j_l}},\dots ,z_{j_{q+1}})}
 \frac{\partial (f_1,\dots ,f_q)}{\partial (z_{j_l},z_{i_1},\dots ,\hat{z_{i_k}},\dots ,z_{i_q})}=0,
\]
which is obvious by expanding the second determinant by the first column and summing over $l$.
\end{proof}

\begin{proposition}
Let $DF(0)\neq 0$. Then $I=J$.
\end{proposition}

\begin{proof}
For any $j_1,\dots ,j_{q+1}$ we get from Lemma \ref{eq1} $DFm_{j_1,\dots ,j_{q+1}}\in I$ and therefore $m_{j_1,\dots ,j_{q+1}}\in I$,
since $DF$ is a unit.
\end{proof}

Now from standard theory we get

\begin{corollary}
\label{cor1}
Let $DF(0)\neq 0$. Then\\
(i) $f_1,\dots ,f_q,\tilde{m}_1,\dots ,\tilde{m}_{n-q}$ defines a regular sequence and
\[
\dim_{\CC}\calA=\RCn
\begin{bmatrix}
\frac{\partial (f_1,\dots ,f_q,\tilde{m}_1,\dots ,\tilde{m}_{n-q})}{\partial (z_1,\dots ,z_n)}\\
f_1\dots f_q \tilde{m}_1 \dots  \tilde{m}_{n-q}
\end{bmatrix}.
\]
(ii) $\RCn
\begin{bmatrix}
\cdot \\
f_1\dots f_q \tilde{m}_1 \dots  \tilde{m}_{n-q}
\end{bmatrix}\colon\calA\to\CC$ defines a linear form.\\
(iii) The induced pairing on $\calA$ is non-degenerate.
\end{corollary}

\subsection{The index of a holomorphic 1-form}
\label{sec22}
From now on we restrict to the case $q=n-2$ and we introduce some more notations. Set $m_i:=m_{1\dots \hat{i}\dots n}$
and for $l\neq k$
\[
f_{l,k}:=\frac{\partial (f_1,\dots ,f_{n-2})}{\partial (z_1,\dots ,\hat{z_l},\dots ,\hat{z_k},\dots ,z_n)}.
\]
Define $M$ to be the matrix
\[
M:=\frac{\partial ((-1)m_1,(-1)^2m_2,\dots ,(-1)^nm_n)}{\partial (z_1,\dots ,z_n)}.
\]
We define $\sigma :=\sigma_2(M)$ as a coefficient of the characteristical polynomial of $M$:
\[
\sum_{i=0}^n(-1)^i\sigma_{n-i}(M)t^i:=\det (M-tI).
\]
We now want to prove a refinement of Corollary \ref{cor1} that will be used to prove our main result.

\begin{lemma}
\label{prop2}
Let $(V,0)$ be smooth and assume $m_1,m_2$ to be a regular $\calO_{V,0}$-sequence. Then
\[
\dim_{\CC}\calA=\RCn
\begin{bmatrix}
DF\cdot \sigma\\ f_1 \dots f_{n-2} m_1 m_2
\end{bmatrix}.
\]
\end{lemma}
Note that we have not assumed $DF(0)\neq 0$. In order to prove this we need a computation.

\begin{lemma}
\label{lem2}
For $j<k$ we have
\[
\frac{\partial (f_1,\dots ,f_{n-2},m_j,m_k)}{\partial (z_1,\dots ,z_n)} =
-f_{j,k}\cdot\sigma \;\mo J.
\]
\end{lemma}

\begin{proof}
The expansion formula gives us
\[
\frac{\partial (f_1,\dots ,f_{n-2},m_j,m_k)}{\partial (z_1,\dots ,z_n)}
= \sum_{1\leq m<l\leq n}(-1)^{m+l+1}f_{m,l}
\frac{\partial (m_j,m_k )}{\partial (z_m,z_l )}.
\]
Now it sufficies to prove $\mo J$
\[
f_{m,l}
\frac{\partial (m_j,m_k )}{\partial (z_m,z_l )} =
f_{j,k}
\frac{\partial (m_m,m_l )}{\partial (z_m,z_l )}.
\]
The first case is $m,l<j$. By Lemma \ref{eq1} we get $\mo J$
\begin{eqnarray*}
f_{m,l}\frac{\partial (m_j,m_k )}{\partial (z_m,z_l )} &=&
f_{m,j}\frac{\partial (m_l,m_k )}{\partial (z_m,z_l )} -
f_{l,j}\frac{\partial (m_m,m_k )}{\partial (z_m,z_l )}\\
&=& f_{m,k} \frac{\partial (m_l,m_j )}{\partial (z_m,z_l )}-
f_{j,k} \frac{\partial (m_l,m_m)}{\partial (z_m,z_l )}\\
&-&  f_{l,k} \frac{\partial (m_m,m_j )}{\partial (z_m,z_l )}+
f_{j,k} \frac{\partial (m_m,m_l )}{\partial (z_m,z_l )}.
\end{eqnarray*}
On the other hand we get again by Lemma \ref{eq1} $\mo J$
\[
f_{m,l}\frac{\partial (m_j,m_k )}{\partial (z_m,z_l )}
= f_{m,k}
\frac{\partial (m_j,m_l )}{\partial (z_m,z_l )}
  -f_{l,k}\frac{\partial (m_j,m_m)}{\partial (z_m,z_l )}.
\]
Summing the two equations we get the Lemma in this case. All other cases can be proved similarly.
\end{proof}

\begin{proof}[Proof of Lemma \ref{prop2}]
Let $k<l$ with $f_{k,l}(0)\neq 0$. By a trivial generalization of Corollary \ref{cor1} we get
\[
\dim_{\CC}\calA=\RCn
\begin{bmatrix}
\frac{\partial (f_1, \dots ,f_{n-2},m_k,m_l)}
{\partial (z_1,\dots ,z_n)}\\
f_1\dots f_{n-2}m_km_l
\end{bmatrix}.
\]
The first and most complicated case is $3\leq k$. Here we get by Lemma \ref{eq1}
\begin{eqnarray*}
f_{k,l}m_1 &=& -f_{1,k}m_l + f_{1,l}m_k\\
f_{k,l}m_2 &=& -f_{2,k}m_l + f_{2,l}m_k.
\end{eqnarray*}
Furthermore it is not hard to compute
\[
f_{2,k}f_{1,l}- f_{1,k}f_{2,l} =-f_{k,l}DF
\]
and therefore we get by the transformation formula for residues and Lemma \ref{lem2}
\begin{eqnarray*}
\dim_{\CC}\calA
&=& \RCn
\begin{bmatrix}
\frac{-DF}{f_{l,k}}\frac{\partial (f_1, \dots ,f_{n-2},m_k,m_l)}
{\partial (z_1,\dots ,z_n)}\\
f_1\dots f_{n-2} m_1 m_2
\end{bmatrix}\\
&=& \RCn
\begin{bmatrix}
DF\cdot\sigma\\
f_1\dots f_{n-2} m_1 m_2
\end{bmatrix}
\end{eqnarray*}
using that the product of $DF$ and any minor is contained in $I$. All other cases can be proved similarly.
\end{proof}

Now we generalize to the singular case.

\begin{lemma}
\label{lem4}
There exists a linear change of coordinates so that $m_1, m_2$ is a regular $\calO_{V,0}$-sequence.
\end{lemma}

\begin{proof}
Let $\phi\colon(\CC^n,0)\to (\CC^n,0)$ be biholomorphic with $\phi (y)=z$ and set $\psi:=\phi^{-1}$. Denote
by $m_i^y$ the minors computed in $y$-coordinates. By standard computations we get for
$i=1,\dots ,n$
\begin{equation}
\label{eq2}
\begin{split}
m_i^y &= \sum_{j=1}^n
\frac{\partial (\phi_1,\dots ,\hat{\phi}_j,\dots ,\phi_n)}{\partial (y_1,\dots ,\hat{y}_i,\dots ,y_n)}
m_j\circ\phi\\
&= \sum_{j=1}^n(-1)^{i+j}\det D\phi\frac{\partial\psi_i}{\partial z_j}\circ\phi\cdot m_j\circ\phi .
\end{split}
\end{equation}
Since $\calO_{V,0}$ is a complete intersection there are complex numbers $c_{11},\dots ,c_{1n}$, $c_{21},\dots ,c_{2n}$ so that
$g_1,g_2$ is a regular $\calO_{V,0}$-sequence where $g_1:=\sum_{i=1}^nc_{1i}m_i$ and
$g_2:=\sum_{i=1}^nc_{2i}m_i$ and therefore $\phi^*(f_1),\dots ,\phi^*(f_{n-2}),\phi^*(g_1),\phi^*(g_2)$ is a regular
$\calO_{\CC^n,0}$-sequence.
Since the vectors $(c_{11},\dots ,c_{1n})$ and $(c_{21},\dots ,c_{2n})$ are linearly independent
we can extend them to a squared matrix of complex numbers $C$ with $\det C=1$. Define $C'$ to be the matrix with
entries $c_{ij}':=(-1)^{i+j}c_{ij}$. Then $\det C'=1$ and therefore $C'$ defines a biholomorphic map. By equation \ref{eq2} it sufficies now to set
$\phi:=(C')^{-1}$ since we have $\phi^*(g_1)=m_1^y$ resp. $\phi^*(g_2)=m_2^y$.
\end{proof}

We call a system of coordinates as in the Lemma a good system of coordinates. Further we define for good coordinates
\[
\RVn \begin{bmatrix} h\\ m_1 m_2 \end{bmatrix} := \RVn \begin{bmatrix} hdz_1\wedge dz_2 \\
m_1 m_2 \end{bmatrix} .
\]
Now we can easily prove our main result.

\begin{proof}[Proof of Theorem \ref{thm1}]
Since the product of $DF$ and any minor is contained in $I$ the residue defines obviously a linear form on $\calA$ by Lemma
\ref{res}. On the other hand we have
\begin{eqnarray*}
\RVn \begin{bmatrix}\sigma\\ m_1 m_2\end{bmatrix}&=& \RCn
\begin{bmatrix}
DF\sigma\\ f_1\dots f_{n-2} m_1 m_2
\end{bmatrix}\\
&=&\lim_{\epsilon\to 0}\sum_i\res_{\CC^n,p_i}
\begin{bmatrix}
DF\sigma\\ f_{\epsilon ,1}\dots f_{\epsilon ,n-2} m_1 m_2
\end{bmatrix}
\end{eqnarray*}
where $f_{\epsilon}:=(f_1-\epsilon_1,\dots ,f_{n-2}-\epsilon_{n-2})$ and we sum over the zeros of $(m_1,m_2)$
on the Milnor fibre $V_{\epsilon}$. We may first ask if $(m_1,m_2)$ can have a zero $p_i$ on $V_{\epsilon}$
when $\ind_{V_{\epsilon,p_i}}\omega = 0$. In this case let $f_{l,k}(p_i)\neq 0$ and we find
\[
\frac{\partial (f_{\epsilon ,1},\dots ,f_{\epsilon ,n-2},m_l,m_k )}{\partial (z_1,\dots ,z_n)}
\in \calO_{\CC^n,p_i}( f_{\epsilon ,1},\dots ,f_{\epsilon ,n-2},m_l ,m_k)
\]
since one of the minors $m_l, m_k$ does not vanish in $p_i$. Cramers rule and the same matrix transformation as in the proof
of Lemma \ref{prop2} show that $DF\sigma \in \calO_{\CC^n,p_i}(f_{\epsilon ,1},\dots ,f_{\epsilon ,n-2},m_1, m_2)$. This means that the residue vanishes
in such a point $p_i$. Therefore the above sum of residues is the sum of the indices of $\omega$ on the Milnor fibre by Lemma
\ref{prop2} which equals $\dim_{\CC}\calA$ by the theorem of Ebeling and Gusein-Zade.
\end{proof}

We cannot expect that the induced pairing is non-degenerate because $\calA$ in general has no $1$-dimensional socle and therefore
there cannot exist a non-degenerate pairing induced by a linear form. We go into more detail in the next subsection.

\subsection{Properties of the residue form}
\label{sec23}
We want to show that the rank of the induced pairing $\beta$ on $\calA$ doesn't depend on the choice of a good system
of coordinates. Let $\phi\colon (\CC^n,0)\to (\CC^n,0)$ be biholomorphic, $\psi :=\phi^{-1}$ and $\phi (y)=z$. We denote
by $m_i^y$ the minors computed in $y$-coordinates and similarly
\[
DF^y:=\frac{\partial (f_1\circ\phi ,\dots ,f_{n-2}\circ \phi )}{\partial (y_3,\dots ,y_n)}.
\]

\begin{lemma}
\label{ann}
Let $(m_1,m_2)$ and $(m_1^y(\psi ),m_2^y(\psi ))$ be regular $\calO_{V,0}$-sequences. Then for any $h\in\OCn$
\[
\RCn\begin{bmatrix} hDF\\ f_1\dots f_{n-2}m_1 m_2\end{bmatrix}=
\RCn\begin{bmatrix} \det (D\phi )(\psi )hDF^y(\psi )\\ f_1\dots f_{n-2}m_1^y(\psi )m_2^y(\psi )\end{bmatrix}.
\]
\end{lemma}

\begin{proof}
The first step is to show the equality if $(V,0)$ is smooth. Let $k<l$ with $f_{k,l}(0)\neq 0$. In the same way as in the
proof of Lemma \ref{prop2} we get
\[
\RCn\begin{bmatrix} hDF\\ f_1\dots f_{n-2}m_1 m_2\end{bmatrix}=
\RCn\begin{bmatrix} hf_{k,l}\\ f_1\dots f_{n-2}m_k m_l\end{bmatrix}.
\]
Lemma \ref{eq1} shows now
\[
f_{k,l}m_j=\left\{
\begin{array}{rl}
-f_{j,k}m_l+f_{j,l}m_k &\mbox{for }j<k\\
f_{j,k}m_l+f_{j,l}m_k &\mbox{for }k<j<l\\
f_{j,k}m_l-f_{j,l}m_k &\mbox{for }l<j
\end{array}\right .
\]
Using this and equation \ref{eq2} we get
\begin{eqnarray*}
f_{k,l}m_1^y(\psi ) &=&\det (D\phi )(\psi )m_k\Bigl(
\sum_{j=1}^{l-1}(-1)^{j+1}\frac{\partial\psi_1}{\partial z_j}f_{j,l}+ \sum_{j=l+1}^{n}(-1)^{j}\frac{\partial\psi_1}{\partial z_j}f_{j,l}\Bigr)\\
& +&  \det (D\phi )(\psi )m_l\Bigl(
\sum_{j=1}^{k-1}(-1)^{j}\frac{\partial\psi_1}{\partial z_j}f_{j,k}+ \sum_{j=k+1}^{n}(-1)^{j+1}\frac{\partial\psi_1}{\partial z_j}f_{j,k}\Bigr)\\
f_{k,l}m_2^y(\psi ) &=&\det (D\phi )(\psi )m_k\Bigl(
\sum_{j=1}^{l-1}(-1)^{j}\frac{\partial\psi_2}{\partial z_j}f_{j,l}+ \sum_{j=l+1}^{n}(-1)^{j+1}\frac{\partial\psi_2}{\partial z_j}f_{j,l}\Bigr)\\
& +&  \det (D\phi )(\psi )m_l\Bigl(
\sum_{j=1}^{k-1}(-1)^{j+1}\frac{\partial\psi_2}{\partial z_j}f_{j,k}+ \sum_{j=k+1}^{n}(-1)^{j}\frac{\partial\psi_2}{\partial z_j}f_{j,k}\Bigr).
\end{eqnarray*}
Therefore we have a matrix $A$ with
\[
\left(\begin{array}{c} f_{k,l}m_1^y(\psi )\\ f_{k,l}m_2^y(\psi )\end{array}\right) = A
\left(\begin{array}{c} m_k\\ m_l \end{array} \right).
\]
Using the formula
\[
f_{i,l}f_{j,k}=\pm f_{k,l}f_{i,j}\pm f_{j,l}f_{i,k}
\]
where $i,l,k,k$ are pairwise disjoint and the signs depend on the position of the indices, it is not hard to compute
\begin{eqnarray*}
\det A &=&\det (D\phi)^2(\psi )f_{k,l}\sum_{1\leq i<j\leq n}(-1)^{i+j+1}\frac{\partial (\psi_1,\psi_2 )}{\partial (z_i,z_j)}f_{i,j}\\
&=& \det(D\phi)(\psi )f_{k,l}DF^y(\psi ).
\end{eqnarray*}
Application of the transformation formula for residues finishes the proof in the smooth case. To generalize to the singular
case we have to show
\small{
\[
\sum_i\res_{\CC^n,p_i}\begin{bmatrix} hDF\\ f_{\epsilon ,1}\dots f_{\epsilon ,n-2}m_1 m_2\end{bmatrix}=
\sum_j\res_{\CC^n,q_j}\begin{bmatrix} \det (D\phi )(\psi )hDF^y(\psi )\\ f_{\epsilon ,1}\dots f_{\epsilon ,n-2}m_1^y(\psi )m_2^y(\psi )\end{bmatrix},
\]
}
\normalsize{
where the first sum is taken over the zeros of $g_1:=(m_1,m_2)$ on the Milnor fibre $V_{\epsilon}$ and the second over
the zeros of $g_2:=(m_1^y(\psi ),m_2^y(\psi ))$ on $V_{\epsilon}$. The residues equal each other at those points which
are common zeros of $g_1$ and $g_2$. Let $p_i$ be a zero of $g_1$ with $g_2(p_i)\neq 0$. Then $\ind_{V_{\epsilon},p_i}\omega =0$.
If $f_{l,k}(p_i)\neq 0$ we have $hf_{k,l}\in\calO_{\CC^n,p_i}(f_{\epsilon ,1},\dots ,f_{\epsilon ,n-2},m_1, m_2)$. The same transformation as
at the beginning of the proof and application of Cramers rule show
$hDF\in\calO_{\CC^n,p_i}(f_{\epsilon ,1},\dots ,f_{\epsilon ,n-2},m_1, m_2)$. Therefore the residue vanishes at such a point.
Similarly one arguments for a zero $q_j$ of $g_2$ with $g_1(q_j)\neq 0$. Since the above sums of residues are equal the limit procedure gives us the equality
of the residues at the origin.

}\end{proof}

We assume that the chosen coordinates as well as the $y$-coordinates are good. We define
\[
\calB :=\frac{\calO_{V,0}}{(m_1,m_2)},\;\; \calC :=\frac{\calB}{\ann_{\calB}(DF)}.
\]
Since we have $\calB(m_3,\dots ,m_n)\subset\ann_{\calB}(DF)$ we find that $\dim_{\CC}\calC\leq\dim_{\CC}\calA$ and the
residue induces a linear form on $\calC$, the induced pairing on $\calC$ is non-degenerate, since
\[
\RVn\begin{bmatrix} hg\\ m_1 m_2\end{bmatrix}=0\; \textnormal{ for all }g\in\OCn
\]
implies
\[
\RCn\begin{bmatrix} hDFg\\ f_1\dots f_{n-2} m_1 m_2\end{bmatrix}=0\; \textnormal{ for all }g\in\OCn .
\]
By duality we get $hDF\in I$ and therefore $h\in\ann_{\calB}(DF)$. This means that $\calC$ has an $1$-dimensional socle
$\soc\calC$. It is generated by the class of $\sigma$: Let $g(0)=0$. Then we get
\begin{eqnarray*}
\RVn\begin{bmatrix} g\sigma\\ m_1 m_2\end{bmatrix} &=&\lim_{\epsilon\to 0}
\sum_i\RCpi\begin{bmatrix} gDF\sigma\\ f_{\epsilon ,1}\dots f_{\epsilon ,n-2} m_1 m_2\end{bmatrix}\\
&=& \lim_{\epsilon\to 0}\sum_ig(p_i)\ind_{V_{\epsilon},p_i}\omega\\
&=& 0
\end{eqnarray*}
and therefore $\sigma\in\soc\calC$. Furthermore we find
$\rg \beta =\dim_{\CC}\calC$
and we can conclude
\[
\dim_{\CC}\soc\calA\leq\dim_{\CC}\calA -\dim_{\CC}\calC +1.
\]
We want to show that $\dim_{\CC}\calC$ does not depend on the choice of a good system of coordinates. Since we have an exact
sequence
\[
0\rightarrow\ann_{\calB}(DF)\rightarrow\calB\stackrel{\cdot DF}{\rightarrow}\calB\rightarrow
\frac{\calB}{\calB (DF)}\rightarrow 0
\]
we get $\dim_{\CC}\calC =\dim_{\CC}\calB (DF)$. Let
\[
\calB':=\frac{\OCn}{(f_1(\phi ),\dots ,f_{n-2}(\phi), m_1^y, m_2^y)}.
\]
We have to show $\dim_{\CC}\calB (DF)=\dim_{\CC}\calB'(DF^y)$. Since we have
\[
\dim_{\CC}\calB'(DF^y)=
\dim_{\CC}\calB''(DF^y(\psi ))
\]
with
\[
\calB'':=\frac{\OCn}{(f_1,\dots ,f_{n-2}, m_1^y(\psi ),m_2^y(\psi ))}
\]
it remains to construct an isomorphism of vector spaces
\[
\varphi\colon\calB (DF)\to\calB''(DF^y\circ\psi ).
\]
We set $\varphi (gDF):=gDF^y(\psi )$. Now Lemma \ref{ann} and duality show that
\[
g\in\ann_{\calB}(DF)\Longleftrightarrow g\in\ann_{\calB''}(DF^y(\psi ))
\]
which means that $\varphi$ is well defined and injective and of course surjective. We want to summarize:

\begin{proposition}
\label{prop}
(i) The residue induces a non-degenerate pairing on $\calC$.\\
(ii) The $1$-dimensional socle of $\calC$ is generated by $\sigma$.\\
(iii) The dimension of $\calC$ does not depend on the choice of a good system of coordinates.\\
(iv) $\rg\beta =\dim_{\CC}\calC$, in particular if $V$ is smooth then $\beta$ is non-degenerate
and $\sigma$ generates the $1$-dimensional socle of $\calA$.\\
(v) $\dim_{\CC}\soc \calA\leq\dim_{\CC}\calA -\dim_{\CC}\calC +1$.
\end{proposition}

\section{Remarks}
\label{sec4}
\subsection{The case of curve singularities}
\label{sec41}
We have only one minor $m$ here and therefore the algebra $\calA$ is a complete intersection. Residue theory can be applied directly
so that we obtain a non-degenerate pairing on $\calA$. But the dimension of $\calA$ can also be expressed by a relative integral
\[
\dim_{\CC}\calA=\frac{1}{2\pi i}\int_{\Sigma}\frac{dm}{m}
\]
by Corollary \ref{cormult}.

\subsection{The general case}
\label{sec42}
We may ask how to generalize our results when we consider a general icis. The problem is that in general there is no choice of coordinates
so that $(\tilde{m}_1,\dots ,\tilde{m}_{n-q})$ is a regular $\calO_{V,0}$-sequence and therefore the residues don't exist. There is
a simple reason that was pointed out to the author by S.M. Gusein-Zade. Let $n>2q+1$. Assume that $f_1,\dots ,f_q,\tilde{m}_1,
\dots ,\tilde{m}_{n-q}$ define an isolated zero. Consider the $(q+1)$ maximal minors of the matrix
\[
\left(\begin{array}{ccc}
\frac{\partial f_1}{\partial z_{n-q+1}} & \dots & \frac{\partial f_1}{\partial z_n}\\
\vdots & \ddots & \vdots\\
\frac{\partial f_q}{\partial z_{n-q+1}} & \dots & \frac{\partial f_q}{\partial z_n}\\
\omega_{n-q+1} & \dots & \omega_n
\end{array}\right) .
\]
Then the vanishing of these minors at the origin would imply that $2q+1$ equations define an isolated zero which is not possible.
If one of these minors is a unit in $\OCn$, then $\mathcal{A}$ is a complete intersection.

\subsection{Linear forms on spaces of relative holomorphic forms}
\label{sec43}
It also seems to be natural to define linear forms on the vector space $\Omega^{n-q}_{V,0}/\omega\wedge\Omega^{n-q-1}_{V,0}$
since it's dimension is also $\IVn\omega$: We have an exact sequence
\[
0\longrightarrow T\Omega^{n-q}_{V,0}\longrightarrow\frac{\Omega^{n-q}_{V,0}}{\omega\wedge\Omega^{n-q-1}_{V,0}}
\stackrel{\lambda}{\longrightarrow}\frac{J'}{J}\longrightarrow 0
\]
and $\dim_{\CC}T\Omega^{n-q}_{V,0}=\dim_{\CC}\OCn /J'$ where $T\Omega^{n-q}_{V,0}$ is the torsion submodule and $J'$
the ideal in $\OCn$ generated by the components of $f$ and the maximal minors of it's Jacobi matrix. Details of this argumentation can
be found in \cite{gre}. If we have a regular $\calO_{V,0}$-sequence $(g_1,\dots ,g_{n-q})$ and the residue
$\RVn\Bigl[\begin{array}{c} \cdot\\ g_1\dots g_{n-q}\end{array}\Bigr]$ vanishes on $\omega\wedge\Omega^{n-q-1}_{V,0}$ the
absolute residue $\RCn\Bigl[\begin{array}{c} \cdot\\ g_1\dots g_{n-q}f_1\dots f_q\end{array}\Bigr]$ vanishes on $J$ since
$m_{j_1,\dots ,j_{q+1}}\in\lambda ( \omega\wedge\Omega^{n-q-1}_{V,0})$ for each minor.

\end{document}